\documentclass[12pt,a4paper]{amsart}
\usepackage{hyperref,amsmath,amsfonts,mathrsfs}
\usepackage{graphics}
\theoremstyle{theorem}

\theoremstyle{definition}

\numberwithin{equation}{section} \makeatletter
\@namedef{subjclassname@2010}{ 2010 Mathematics Subject
Classification} \makeatother
\begin{document}
\title[energy of signed digraphs]{ENERGY OF SIGNED DIGRAPHS}
\author{S. Pirzada and Mushtaq A. Bhat}

\address{Department of Mathematics \\
    University of Kashmir \\
   Srinagar, Hazratbal 190006\
   \\ India}
 \email{pirzadasd@kashmiruniversity.ac.in; sdpirzada@yahoo.co.in}
 \email{mushtaqab1125@gmail.com}
\date{}

\subjclass[2000]{fill it,}
\begin{abstract}
In this paper we extend the concept of energy to signed digraphs. We obtain Coulson's integral formula for energy of signed digraphs. Formulae for energies of signed directed cycles are computed and it is shown that energy of non cycle balanced signed directed cycles increases monotonically with respect to number of vertices. Characterization of  signed digraphs having energy equal to zero is given. We extend the concept of non complete extended $p$ sum (or briefly, NEPS) to signed digraphs. An infinite family of equienergetic signed digraphs is constructed. Moreover, we extend McClelland's inequality to signed digraphs and also obtain sharp upper bound for energy of signed digraph in terms the number of arcs. Some open problems are also given at the end.
\end{abstract}
\maketitle{}
\section{introduction}\label{sec1}
A signed digraph (or briefly sidigraph) is defined to be a pair $S=(D,\sigma)$ where $D=(V,\mathscr{A})$ is the underlying digraph and $\sigma: \mathscr{A}\rightarrow \{-1, 1\}$ is the signing function. The sets of positive and negative arcs of $S$ are respectively denoted by $\mathscr{A}^+(S)$ and $\mathscr{A}^-(S)$. Thus $\mathscr{A}(S)=\mathscr{A}^+(S)\cup \mathscr{A}^- (S)$. A sidigraph is said to be homogeneous if all of its arcs have either positive sign or negative sign, otherwise heterogeneous. A sidigraph without orientations of its arcs (when it is regarded as a simple undirected sigraph) is known as its underlying sigraph and is denoted by $S^u$.  \\
\indent Two vertices are adjacent if they are connected by an arc. A path of length $n-1$ $(n\geq2)$, denoted by $P_n$, is a sidigraph  on $n$ vertices $v_1, v_2,\cdots,v_n$ with $n-1$  signed arcs $(v_i,v_i+1)$. A cycle of length $n$, is a sidigraph having vertices $v_1, v_2,\cdots,v_n$ and signed arcs $(v_i, v_i+1)$ $i=1,2,\cdots,n-1$ and $(v_n,v_1)$. A sidigraph is linear if each vertex has both indegree and outdegree equal to one. The sign of a sidigraph is defined as the product of signs of its arcs. A sidigraph is said to be positive (negative) if its sign is positive (negative) i.e., it contains an even (odd) number of negative arcs. A signed digraph is said to be all-positive (respectively, all negative) if all its arcs are positive (negative). A signed digraph is said to be cycle balanced if each of its cycles is positive, otherwise non cycle balanced. Throughout this paper we call cycle balanced cycle a positive cycle and non cycle balanced cycle a negative cycle and respectively denote them by $C_n$ and ${\bf C}_n$, where $n$ is number of vertices.\\
 A sidigraph is symmetric if $(u,v)\in \mathscr{A}^{+} (S) ~~(or~~\mathscr{A}^-(S))$ then $(v,u)\in \mathscr{A}^{+} (S)~~(or~~\mathscr{A}^-(S))$, where $u,v\in V(S)$. A one to one correspondence between sigraphs and symmetric sidigraphs is given by $S\rightsquigarrow \overleftrightarrow{S}$, where $\overleftrightarrow{S}$ has the same vertex set as that of sigraph $S$, and each signed edge $(u,v)$ is replaced by a pair of symmetric arcs $(u,v)$ and $(v,u)$ both with same sign as that of edge $(u,v)$. Under this correspondence a sigraph can be identified with a symmetric sidigraph. A sidigraph is said to be skew symmetric if its adjacency matrix is skew symmetric. A skew symmetric sidigraph is obtained by replacing each edge of sigraph by a pair of opposite oriented and opposite signed arcs. We denote a skew symmetic digraph of order $n$ by ${\bf S}_n$.  \\
\indent  Two sidigraphs (sigraphs) on same number of vertices are said to be cospectral if they have same spectrum. Acharya, Gill, Patwardhan \cite{agp} introduced the concept of Quasicospectrality in digraphs (graphs). Two digraphs (graphs) $G$ and $H$ of same order are said to be quasicospectral if there exist sidigraphs (sigraphs) $S$ and $ S'$ on $G$ and $H$ respectively which are cospectral. The weighted directed graph $S$ of an $n \times n$ matrix $M=(m_{ij})$ of reals  consists of $n$ vertices, with vertex $i$ joined by a directed arc with weight $m_{ij}$ to vertex $j$ if and only if $m_{ij}$ is non-zero. In case the matrix consists of entries $-1,0$ and $1$, then we get a sidigraph. Thus there is a one to one correspondence between the set of integral $(-1,0,1)$-matrices of order $n$ and  set of sidigraphs of order $n$. \\

\indent The cartesian product (or sum) of two sidigraphs $S_1=(V_1,A_1,\sigma_1)$ and  $S_2=(V_2,A_2,\sigma_2)$ denoted by $S_1\times S_2$ is the sidigraph $(V_1\times V_2,A,\sigma)$, where the arc set is that of the cartesian product of unsigned digraphs and the sign function is defined by\\
 $$\sigma((u_i,v_j),(u_k,v_l))=\left \{\begin{array}{lr}\sigma_1(u_i,u_k), &\mbox{if $j=l$},\\
 \sigma_2(v_j,v_l), &\mbox{if $i=k.$}
 \end{array} \right.$$\\
\indent The Kronecker product (strong product or conjunction) of two sidigraphs $S_1=(V_1,A_1,\sigma_1)$ and  $S_2=(V_2,A_2,\sigma_2)$ denoted by $S_1 \otimes S_2$ is the sidigraph $(V_1\times V_2,A,\sigma)$, where arc set is the arc set of underlying unsigned digraphs and the sign function is defined by $\sigma((u_i,v_j),(u_k,v_l))=\sigma_1(u_i,u_k)\sigma_2 (v_j,v_l).$\\

\indent The adjacency matrix of a sidigraph $S$ whose vertices are ${v_1,v_2,\cdots,v_n}$ is the $n\times n$ matrix $A(S)=(a_{ij})$, where\\
$$a_{ij}=\left\{\begin{array}{lr}\sigma(v_i,v_j), &\mbox{if there is an arc from $v_i$ to $v_j,$}\\
0, &\mbox {otherwise.}
\end{array}\right.$$ \\
The characteristic polynomial $|xI-A(S)|$ of the adjacency matrix $A(S)$  of sidigraph $S$ is called the characteristic polynomial of $S$ and is denoted by $\phi_S(x)$. The eigenvalues of $A(S)$ are called eigenvalues of $S$. The following is the coefficient Theorem for sidigraphs \cite{agp}\\

\noindent {\bf Theorem 1.1.} If $S$ be a sidigraph with characteristic polynomial 
$$\phi_S(x)=x^n+c_1 x^{n-1} +\cdots+c_{n-1}x+c_n$$\\
then $$ c_i=\sum\limits_{L\in \pounds_i}(-1)^{p(L)}\prod\limits_{Z\in c(L)} s(Z),$$
for all $i=1,2,\cdots,n$, where $\pounds_i$ is the set of all linear subdigraphs  $L$ of $S$ of order $i$, $p(L)$ denotes number of components of $L$ and $c(L)$ denotes the set of all cycles of $L$ and $s(Z)$ the sign of cycle $S$.\\

\indent If $S$ is an undirected sigraph then $S$ can be viewed as a digraph $\bar{S}$ by identifying each signed edge of $S$ with a directed cycle of length $2$ with each arc having same sign as that of the corresponding edge. For undirected sigraphs Theorem $1.1$ takes the following form \cite{ga}.\\

\noindent {\bf Theorem 1.2.} If $S$ is a sigraph with characteristic polynomial
$$\phi_S(x)=x^n+c_1 x^{n-1} +\cdots+c_{n-1}x+c_n $$
then $$ c_i=\sum\limits_{L\in \pounds_i}(-1)^{p(L)} 2^{c(L)}\prod\limits_{Z\in c(L)} s(Z),$$\\
for all $i=1,2,\cdots,n$, where $\pounds_i $ is the set of all basic figures $L$ of $S$ of order $i$, $p(L)$ denotes number of components of $L$ and $c(L)$ the number of cyclic components of $L$.\\

\indent The spectral criterion for cycle balance of sidigraphs (sigraphs) given by Acharya is as follows.\\

\noindent {\bf Theorem 1.3}\cite{a,agp}. A sidigraph (sigraph) is cycle balanced (balanced) if and only if it is cospectral with the underlying unsigned digraph (graph).\\

\section{Energy of sidigraphs}
For spectra and energy of graphs and digraphs see \cite{br,cds,ig1,ig2,pr}. Germina, Hameed and Zaslavsky \cite{gh1} defined energy of a sigraph to be the sum of absolute values of sigraph eigenvalues. In this section, we extend the concept of energy to sidigraphs. Unlike sigraphs the adjacency matrix of a sidigraph need not be real symmetric, so eigenvalues can be complex numbers.\\

\noindent {\bf Definition 2.1.}
Let $S$ be a sidigraph of order $n$ having eigenvalues $z_1,z_2,\cdots,z_n$. The energy of $S$ is defined as
$$E(S)=\sum\limits_{j=1}^{n}|\Re z_j|,$$
where  $\Re z_j $ denotes real part of complex number $z_j. $\\

\indent If $S$ is a signed graph and $\overleftrightarrow{S}$ be its symmetric sidigraph, then clearly  $A(S)=A(\overleftrightarrow{S})$ and so $E(S)=E(\overleftrightarrow{S})$. In this way, definition $2.1$ generalizes the concept of energy of  undirected sigraphs.\\

\noindent {\bf Example 2.2.}
Let $S$ be a sidigraph shown in Figure $1$. Dotted arrows denote negative arcs and plane ones denote positive arcs. Clearly, $S$ is non cycle balanced sidigraph. By Theorem $1.1$, the characteristic polynomial of $S$ is $ \phi_S(x)=x^{10}+x^7=x^7(x^3+1)$. The eigenvalues of $S$ are $0^7,1,\frac{1\pm \sqrt{3}\iota}{2}$, where $\iota=\sqrt{-1}$, so $E(S)=2$.\\

\unitlength 1mm 
\linethickness{0.4pt}
\ifx\plotpoint\undefined\newsavebox{\plotpoint}\fi 
\begin{picture}(66,56.75)(0,0)
\put(30.5,17.5){\circle*{1}}
\put(46.75,17.5){\circle*{1}}
\put(30.75,32.75){\circle*{1.12}}
\put(46.75,32.5){\circle*{1}}
\put(47,51.5){\circle*{1}}
\put(65.25,32.25){\circle*{1.12}}
\put(65.25,51.5){\circle*{1}}
\put(65.75,51.5){\circle*{.5}}
\put(46.5,25.25){\vector(0,1){.07}}\put(46.5,17.75){\line(0,1){15}}
\put(46.75,42.38){\vector(0,-1){.07}}\put(46.75,51.75){\line(0,-1){18.75}}
\put(55.75,32.75){\vector(1,0){.07}}\put(46.75,32.75){\line(1,0){18}}
\put(65.25,42.25){\vector(0,1){.07}}\multiput(65.18,32.18)(0,.9524){22}{{\rule{.4pt}{.4pt}}}
\put(56,52){\vector(-1,0){.07}}\multiput(65.18,52.18)(-.9737,-.0263){20}{{\rule{.4pt}{.4pt}}}
\put(38.5,32.88){\vector(-1,0){.07}}\multiput(46.43,32.68)(-.9412,.0147){18}{{\rule{.4pt}{.4pt}}}
\put(30.38,25.38){\vector(0,-1){.07}}\multiput(30.43,32.93)(-.0156,-.9531){17}{{\rule{.4pt}{.4pt}}}
\put(38.25,17.63){\vector(1,0){.07}}\multiput(30.18,17.68)(.9412,-.0147){18}{{\rule{.4pt}{.4pt}}}
\put(39.88,39.5){\vector(-1,1){.07}}\multiput(46.43,32.68)(-.631,.6429){22}{{\rule{.4pt}{.4pt}}}
\put(33.5,46.25){\circle*{1}}
\put(28.5,56.25){\circle*{1}}
\put(21.5,47.25){\circle*{1}}
\put(24.75,52){\vector(3,4){.07}}\multiput(21.18,47.43)(.5385,.6923){14}{{\rule{.4pt}{.4pt}}}
\put(27.63,46.63){\vector(-1,0){.07}}\multiput(33.75,46.5)(-1.53125,.03125){8}{\line(-1,0){1.53125}}
\put(30.88,51.75){\vector(2,-3){.07}}\multiput(28.25,56.25)(.03365385,-.05769231){156}{\line(0,-1){.05769231}}
\put(42,12){\makebox(0,0)[cc]{S}}
\put(36.75,5){\makebox(0,0)[cc]{Figure 1}}
\end{picture}

\noindent{\bf Example 2.3.}
Let $S$ be an acyclic sidigraph. Then by Theorem $1.1$, the characteristic polynomial of $S$ is $\phi(S)=x^n$, so that $spec(S)=\{0^n\}$ and hence $E(S)=0$.\\

\noindent{\bf Example 2.4.} Consider ${\bf S}_n$, the skew symmetric sidigraph on $n\ge 2$ vertices, then by Theorem $1.1$, $\phi_{{\bf S}_n}(x)=x^{n-2}(x^2+(n-1))$ so that eigenvalues are $0^{n-2}, \pm\iota\sqrt{n-1}$. Thus all skew symmmetric sidigraphs are cospectral with energy equal to zero.\\

\noindent {\bf Example 2.5.}
If $S$ is the sidirected cycle on $n$ vertices, then the characteristic polynomial of $S$ is $\phi(S)=x^n+(-1)^{[s]}$, where the symbol $[s]$ is defined as $[s]=1$ or $0$ according as $S$ is positive or negative. If $S=C_n$, then eigenvalues are $e^{\frac{2\iota j\pi}{n}},~ j=0,1,\cdots,n-1$. So, $E(S)=\sum\limits_{j=0}^{n-1}|\cos(\frac{2j\pi}{n})|$. If $S={\bf C}_n$, then the eigenvalues are $e^{\frac{\iota (2j+1)\pi}{n}},~ j=0,1,\cdots,n-1$. Therefore, $E(S)=\sum\limits_{j=0}^{n-1}|\cos(\frac{(2j+1)\pi}{n})|$. In particular if $S={\bf C}_4$, the eigenvalues are $\frac {1\pm \iota}{\sqrt{2}}, \frac {-1\pm \iota}{\sqrt{2}}$ and $E(S)=2\sqrt{2}$.\\

\noindent {\bf Example 2.6.}
Let $S$ be a sidigraph having $n$ vertices and unique cycle of length $r$, where $2\leq r\leq n$. Then by Theorem $1.1$, $\phi_S(x)=x^n+(-1)^{[s]} x^r=x^r(x^{n-r}+(-1)^{[s]})$, where the symbol $[s]$ is defined as $[s]=1$ or $0$ according as $S$ is cycle balanced or non cycle balanced. Clearly, energy equals to the energy of the unique cycle.\\

\indent Given $t$ sidigraphs $S_1,S_2,\cdots S_t$, their direct product denoted by $S_1\oplus S_2 \oplus \cdots \oplus S_t $ is the sidigraph with $ V_{S_1 \oplus S_2 \oplus \cdots \oplus S_t}=\bigcup \limits_{j=1}^{t} V(S_j)$ and arc set $\mathscr{A}_{S_1\oplus S_2 \oplus \cdots \oplus S_t}=\bigcup \limits_ {j=1}^{t} \mathscr{A}(S_j)$.\\

\noindent {\bf Theorem 2.7.}
Let $S$ be a sidigraph on n vertices and $S_1,S_2,\cdots,S_k$ be its strong components. Then
$E(S)=\sum\limits_{j=1}^{k}E(S_j).$

\noindent{\bf Proof.}
Let $Y=\{a \in A(S): a\notin c(S)\}$, where $c(S)$ is the set of all cycles of $S$. By Theorem $1.1$, $ \phi_S(x)=\phi_{S-Y}(x)$, where $S-Y$ is the sidigraph obtained from $S$ by deleting the non-cyclic arcs. Clearly, $S-Y=S_1\oplus S_2\oplus\cdots\oplus S_k$ and adjacency matrix of sidigraph $S-Y$ is in block diagonal form with diagonal blocks as the adjacency matrices of strong components. Therefore $\phi_{S-Y}(x)=\phi_{S_1}(x)\phi_{S_2}(x)\cdots\phi_{S_k}(x)$
and so $E(S)=\sum\limits_{j=1}^{k}E(S_j)$. \qed \\

\noindent {\bf Remark 2.8.}
From Theorem $1.1$,  $c_i=\sum\limits_{L\in \pounds_i } (-1)^{p(L)}, i=1,2,\cdots,n$. Clearly, this sum contains positive and negative ones.\\
$+1$ arises if and only if\\ $(a)$ Number of components of $L\in \pounds_i$ are odd  and $s(L)<0$. We call such linear sidigraphs as type $a$ linear sidigraphs.\\
$(b)$ Number of components of $L\in \pounds_i$ is even and $s(L)> 0$. We call such linear sidigraphs as type $b$.\\
$-1$ will occur if and only if\\ $(c)$ Number of components of $L\in \pounds_i$ is odd and $s(L)> 0$. We call such linear sidigraphs as type $c$.\\ $(d)$ Number of components of $L\in \pounds_i$ is even and $s(L)< 0$. We call such linear sidigraphs as type $d$.\\
\indent From above we observe that $c_i=0$ if and only if either $S$ is acyclic or for each $i$, number of linear sidigraphs of order $i$ of type $a$ or type $d$ is equal to the number of linear sidigraphs  of order $i$ of type $b$ or type $c$.\\

\indent An immediate consequence of Remark $2.7$ is the following Lemma.\\

\noindent {\bf Lemma 2.9.} An integral $(-1,0,1)$ matrix is nilpotent if and only if its underlying sidigraph $S$ is either acyclic or in $S$, for each $i=1,2,\cdots,n$, the number of linear sidigraphs of order $i$ of type $a$ or type $d$ is equal to number of linear sidigraphs  of order $i$ of type $b$ or type $c$.  \qed \\

\indent Unlike unsigned strong component, energy of a signed directed strong component can be zero, for example, sidigraph $S_1$ in Figure 2. The following result characterizes sidigraphs having energy equal to zero.\\

\noindent{\bf Theorem 2.10.} Let $S$ be a sidigraph of order $n$. Then $E(S)=0$ if and only if $S$ satisfies one of the following conditions
$(I)$ $S$ is acyclic $(II)$ $S$ is skew symmetric $(III)$  for each $i=1,2,\cdots,n$, the number of linear sidigraphs of order $i$ of type $a$ or type $d$ is equal to number of linear sidigraphs of order $i$ of type $b$ or type $c$.\\
\noindent {\bf Proof.}
Let $S$ be a sidigraph of order $n$. If $S$ is acyclic, then by Lemma $2.9$, $\phi_S (x)=x^n$ and so $E(S)=0.$ If $S$ is skew symmetric sidigraph, then the eigenvalues of $S$ are $0^{n-2}, \pm\iota\sqrt{n-1}$, therefore $E(S)=0$. If $S$ satisfies $(III)$, then by Lemma $2.9$, $\phi_S (x)=x^n$, so that $E(S)=0.$\\
\indent ~~ Conversely, if $E(S)=0$, then either eigenvalues of $S$ are $0^{n-2}, \pm\iota\sqrt{n-1}$ or all equal to zero. Using the fact that a real matrix is skew symmetric if and only if all its eigenvalues are of the form $\pm \iota \alpha$, where $\alpha \in  \Bbb R$ and by Lemma $2.9$, the result follows. \qed \\

\section {Computation of energy of signed directed cycles}

We first calculate energy formulae for positive cycles. Let $C_n$ be a positive cycle on $n\ge 2$ vertices. The characteristic polynomial of $C_n$ is $\phi_{C_n} (x)=x^n-1$, so that the eigenvalues are  $e^{\frac{2\pi \iota j}{n}}, j=0,1,\cdots,n-1,$ $\iota = \sqrt{-1}$. Consequently energy of $C_n$ is\\
$$E(C_n)=\sum\limits_{j=0}^{n-1}|\cos\frac{2j\pi}{n}|.$$
Given a positive integer $n$, it has one of the forms $4k$, or $2k+1$, or $4k+2$, where $k\ge 0$. \\

If $n=4k$, then
\begin{align*}
E(C_n)&=\sum\limits_{j=0}^{4k-1}|\cos\frac{2j\pi}{4k}| =\sum\limits_{j=0}^{4k-1}|\cos\frac{j\pi}{2k}|
=2\sum\limits_{j=0}^{2k-1}|\cos\frac{j\pi}{2k}|\\& =2+4\sum\limits_{j=1}^{k-1}\cos\frac{j\pi}{2k}
=2+4\{ \frac{-1}{2}+\frac{\sin\frac{(k-\frac{1}{2}) \pi}{2k}}{2 \sin \frac{\pi}{4k}}\}=2\cot\frac{\pi}{n}.
\end{align*}

If $n=2k+1$, then
\begin{align*}
E(C_n)&=\sum\limits_{j=0}^{2k}|\cos\frac{2j\pi}{2k+1}|=1+2\sum\limits_{j=1}^{k}|\cos\frac{2j\pi}{2k+1}|
=1+2\sum\limits_{j=1}^{k}\cos\frac{j\pi}{2k+1}\\&=1+2\{ \frac{-1}{2}+\frac{\sin\frac{(k+\frac{1}{2}) \pi}{2k+1}}{2 \sin \frac{\pi}{2(2k+1)}}\}=2\csc\frac{\pi}{2n}.
\end{align*}

If $n=4k+2$ then
\begin{align*}
 E(C_n)&=\sum\limits_{j=0}^{4k+1}|\cos\frac{2j\pi}{4k+2}|=\sum\limits_{j=0}^{4k+1}|\cos\frac{j\pi}{2k+1}|
=2\sum\limits_{j=0}^{2k}|\cos\frac{j\pi}{2k+1}|\\& =2+4\sum\limits_{j=1}^{k}\cos\frac{j\pi}{2k+1}
=2+4\{ \frac{-1}{2}+\frac{\sin\frac{(k+\frac{1}{2}) \pi}{2k+1}}{2 \sin \frac{\pi}{4k+2}}\}=2\csc\frac{\pi}{n}.
\end{align*}

\indent We now calculate exact formulae for the energy of negative cycles of length $n$. Let ${\bf C}_n$ denote the negative cycle with $n$ vertices. Then $\phi_{{\bf C}_n}(x)=x^n+1$ and so $Spec({\bf C}_n)=e^{\frac{(2j+1)\pi \iota}{n}}, j=0,1,\cdots,n-1$, $\iota = \sqrt{-1}$. Therefore energy is given by\\
$$E({\bf C}_n)=\sum\limits_{j=0}^{n-1}|\cos\frac{(2j+1)\pi}{n}|.$$
If $n=4k$, then

\begin{align*}
 E({\bf C}_n)&=\sum\limits_{j=0}^{4k-1}|\cos\frac{(2j+1)\pi}{4k}|=2\sum\limits_{j=0}^{2k-1}|\cos\frac{(2j+1)\pi}{4k}|
=4\sum\limits_{j=0}^{k-1}\cos\frac{(2j+1)\pi}{4k}\\&=4\{ \cos\frac{\pi}{4k}+ \cos\frac{3\pi}{4k}+\cdots + \cos\frac{(2k-1)\pi}{4k}\}=4\frac{\cos(\frac{\pi}{4k}+\frac{k-1}{2} \frac{2\pi}{4k}) \sin k\frac{2\pi}{8k}}{ \sin\frac{2\pi}{8k}}=2\csc \frac {\pi}{n}.
\end{align*}

If $n=4k+2$, then
\begin{align*}
 E({\bf C}_n)&=\sum\limits_{j=0}^{4k+1}|\cos\frac{(2j+1)\pi}{4k+2}|=4\sum\limits_{j=0}^{k-1}\cos\frac{(2j+1)\pi}{4k+2}
\\&=4\{ \cos\frac{\pi}{4k+2}+ \cos\frac{3\pi}{4k+2}+\cdots + \cos\frac{(2k-1)\pi}{4k+2}\}=2\cot \frac {\pi}{n}.
\end{align*}

If $n=2k+1$, then since $-1$ is the eigenvalue of ${\bf C}_n$, we have $spec({\bf C}_n)$ $=-spec(C_n)$, and so $E({\bf C}_n)=E(C_n)$.\\

\indent Summarizing, all the above cases can be written as follows: \\

$$E(C_n)=\left \{\begin{array}{lr} 2\cot\frac {\pi}{n}, &\mbox{if $n=4k$},\\
\csc\frac {\pi}{2n}, ~~~~&\mbox{if $n=4k+1$~or~$n=4k+3$},\\
2\csc\frac {\pi}{n}, &\mbox{if $n=4k+2$}\\
\end{array} \right. $$

and\\

$$E(\bf{C}_n)=\left \{\begin{array}{lr} 2\csc\frac {\pi}{n}, &\mbox{if $n=4k$},\\
\csc\frac {\pi}{2n}, ~~~~&\mbox{if $n=4k+1$~or $n=4k+3$},\\
2\cot\frac {\pi}{n}, &\mbox{if $n=4k+2$}.\\
\end{array} \right. $$

Pena and Rada \cite{pr} proved that the energy of directed unsigned cycles increases monotonically with respect to order $n$. From energy formulae for positive and negative sidirected cycles, the following two results are immediate.\\

\noindent {\bf Theorem 3.1.}
Among all non cycle balanced unicyclic sidigraphs on $n$ vertices, the cycle has the largest energy. Moreover, energy of negative cycles increases monotonically with respect to the order. Minimal energy is attained in ${\bf C}_2 $, a negative cycle of order $2$. \qed \\

\noindent {\bf Theorem 3.2.}
Energy of positive and negative cycles satisfy the following:\\
$(1)$ Energy of positive cycle of odd order equals energy of negative cycle of same order.\\
$(2)$ Energy of negative cycle of even order is greater than energy of positive cycle of same order if and only if $n=4k$.\\
$(3)$ Energy of negative cycle of even order is less than energy of positive cycle of same order if and only if $n=4k+2$.  \\

One of the fundamental results in theory of graph energy is Coulsons integral formula. We now state Coulsons integral formula for energy of sidigraphs. The proof is similar as in digraphs.\\

\noindent {\bf Theorem 3.3.}
Let $S$ be a sidigraph with $n$ vertices having characteristic polynomial $\phi_{S} (x)$. Then
$$ E(S)=\sum\limits_{j=1}^{n}|\Re z_j|=\frac{1}{\pi} \int\limits_{-\infty}^{\infty}(n-\frac{\iota x \phi'_{S}(\iota x)}{\phi_{S}(\iota x)})dx,$$ \\
where $z_1,z_2,\cdots,z_n$ are the eigenvalues of sidigraph $S$ and $\int\limits_{-\infty}^{\infty} F(x) dx$ denotes principle value of the respective integral. \\

\indent The Coulson's formula given above is a motivation to define energy of sidigraphs to be the sum of absolute values of real parts of eigenvalues. \\

\noindent {\bf Example 3.4.}
Consider the cycle ${\bf C}_4 $, the characteristic polynomial is $\phi_{{\bf C}_4}(x)=x^4+1$ and hence\\
$E({\bf C}_4)= \frac{1}{\pi} \int\limits_{-\infty}^{\infty}[4-\frac{4\iota x (\iota x)^3}{(\iota x)^4+1}]dx$
$=\frac{1}{\pi} \int \limits_{-\infty}^{\infty} \frac {4}{x^4+1}dx $  $=\frac{4}{\pi} \frac{\pi}{2\sin \frac{\pi}{4}}=2\sqrt{2}$, as calculated in example $2.5$.  \\

\indent An immediate consequence of Coulson's integral formula is the following observation, the proof of which is similar to the case in unsigned digraphs.\\

\noindent {\bf Corollary 3.5.}
If $S$ is a sidigraph on $n$ vertices, then\\
$$E(S)=\frac{1}{\pi} \int\limits_{-\infty}^{\infty}\frac{1}{x^2} \log(x^n \phi_G(\frac{\iota}{x})dx.$$

Pena and Rada \cite{pr} considered the problem of increasing property of the energy in digraphs, It is natural to consider the same problem for sidigraphs.\\
\indent Let ${S}_{n,h}$ denote the set of sidigraphs with $n$ vertices and every cycle having length $h$. First we shall calculate characteristic polynomial of such sidigraphs.\\

\noindent {\bf Theorem 3.6.}
If $ S\in S_{n,h} $, then $ \phi_S (x)=x^n+\sum\limits_{j=1}^{\lfloor{\frac{n}{h}}\rfloor} (-1)^k c^* (S,kh) x^{n-kh}, $ where $c^* (S,kh) =$ number of  positive linear sidigraphs of order $kh$ $-$ number of  negative linear sidigraphs of order $kh$, for every $k=1,2,\cdots,{\lfloor{\frac{n}{h}}\rfloor} $.\\
\noindent {\bf Proof.} By Theorem $1.1$, the coefficient $c_{kh} = \sum\limits_{L\in\pounds_{kh}} (-1)^{p(L)}\prod\limits_{Z\in c(L)} s(Z) $. As each cycle of $S$ has length $h$, it follows that $c_{kh}=(-1)^k c^*(S,kh)$, where $c^* (S,kh) =$ number of  positive linear sidigraphs consisting of $k$ cycles each of length $h$ $-$ number of  negative linear sidigraphs consisting of $k$ cycles each of length $h$. Clearly, $c_j=0$ if $j$ is not a multiple of $h$, because $\pounds_j =\phi$ in that case. \qed \\

\noindent {\bf Definition 3.7.} Let $S^*_{n,h}=\{S\in S_{n,h}: c^* (S,kh)\ge 0, k=1,2,\cdots,\lfloor \frac{n}{h}\rfloor\}$. Define a quasi order relation over $S^*_{n,h}$ as follows. Let $S_1$ and $S_2$ be two elements of $S^*_{n,h}$. Define $S_1\preceq S_2$ if for all $k=1,2,\cdots,\lfloor \frac{n}{h}\rfloor$, $c^* (S_1,kh)\leq c^* (S_2,kh)$. If $S_1\preceq S_2$ and there exists $k$ such that $c^* (S_1,kh) < c^* (S_2,kh)$ then $S_1\prec S_2$. Clearly, this is a reflexive and transitive relation over  $S^*_{n,h}$.\\

\noindent {\bf Theorem 3.8.} Let $h$ be an integer of the form $4t+2,~t\ge 0$. The energy increases with respect to quasi-order relation defined over $S^*_{n,h}$, i.e., if $S_1,~S_2\in S^*_{n,h}$ then $S_1\prec S_2 \Rightarrow E(S_1) < E(S_2).$\\
\noindent{\bf Proof.}
Let $H\in S^*_{n,h} $, by Theorem $3.6$, we have
 $$\phi_H (x)=x^n+\sum\limits_{j=1}^{\lfloor{\frac{n}{h}}\rfloor} (-1)^k c^* (H,kh) x^{n-kh},$$
 so that $$\phi_H (\frac{\iota}{x})=\frac{\iota^n}{x^n}+\sum\limits_{j=1}^{\lfloor{\frac{n}{h}}\rfloor} (-1)^k c^* (H,kh) \frac {\iota^{n-kh}}{x^{n-kh}}$$\\
 $$=\frac{\iota ^n}{x^n}(1+\sum\limits_{j=1}^{\lfloor{\frac{n}{h}}\rfloor} (-1)^k c^* (H,kh)x^{kh}\iota ^{-k(4t+2)})$$\\
 $$=\frac{\iota ^n}{x^n}(1+\sum\limits_{j=1}^{\lfloor{\frac{n}{h}}\rfloor} c^* (H,kh)x^{kh}).$$\\
 Applying Corollary $3.5$, we have\\
 $E(S)=\frac{1}{\pi} \int\limits_{-\infty}^{\infty}\frac{1}{x^2} \log(x^n\frac{\iota^n}{x^n}(\sum\limits_{j=1}^{\lfloor{\frac{n}{h}}\rfloor} c^* (H,kh)x^{kh}))dx=$ $\frac{1}{\pi} \int\limits_{-\infty}^{\infty}\frac{1}{x^2} \log(\iota^n(1
 +\sum\limits_{j=1}^{\lfloor{\frac{n}{h}}\rfloor} c^* (H,kh)x^{kh}))dx$
 $=\frac{1}{\pi} \int\limits_{-\infty}^{\infty}\frac{1}{x^2} \log(1
  +\sum\limits_{j=1}^{\lfloor{\frac{n}{h}}\rfloor} c^* (H,kh)x^{kh})dx $, as $\frac{1}{\pi} p.v. \int\limits_{-\infty}^{+\infty}\log(\iota^n)\frac{dx}{x^2}=0$, where $p.v.$ stands for principal value of Cauchy's integral. It is clear from above energy expression that energy increases with respect to quasi-order relation defined over $S^*_{n,h}$. \qed \\

\section {NEPS in sidigraphs}

Recall that \cite{hj} Kronecker product of two matrices $A=(a_{ij})_{r\times s}$ and $B=(b_{ij})_t\times u$ denoted by $A\otimes B$ is matrix of order $rt\times su$ obtained by replacing each entry $a_{ij}$ of $A$ by a block $a_{ij}B$. Thus $A\otimes B$ consists of all $rtsu$ possible products of an entry of $A$ with an entry of $B$. The Kronecker product is a componentwise operation, i.e., $(A\otimes B)(C\otimes D)=(AC)\otimes(BD)$, provided the products $AC$ and $BD$ exist. This operation is also associative, so we can define the multiple product $A_1\otimes A_2\otimes\cdots\otimes A_m$. Let order of $A_i$ be $r_i\times s_i$. We index elements of $A_i$ by $a_i;jk$, and those of multiple product by a pair of $m-$tuples, a row index $ j=(j_1,j_2,\cdot,j_m)$ and a column index $ k=(k_1,k_2,\cdots,k_m)$, where $1\le j_i\le r_i$ and $1\le k_i\le s_i$. The element $a_{jk}$ of the product matrix is $a_{jk}=a_{1;j_1 k_1} a_{2;j_2 k_2}\cdots a_{m;j_m k_m}.$\\

\noindent {\bf Lemma 4.1} \cite{gh2}. Let $A_i, i=1,2,\cdots,m$, be a square matrix of order $n_i$ and $\xi_{ij},  ~j=1,2,\cdots,n_i$ be its eigenvalues. Let $k_1,k_2,\cdots, k_m $ be non-negative integers, then the $n_1 n_2\cdots n_m$ eigenvalues of the matrix $A^{k_1}_1\otimes \cdots \otimes A^{k_m}_m$ are $\xi_{j_1 j_2 \cdots j_m}=\xi^{k_1}_{1j_1}\cdots \xi^{k_m}_{1j_m}$ for $1\le j_i \le n_i.$\\

 Let $k_p=(k_{p1},k_{p1},\cdots,k_{pm}), p=1,2,\cdots, q$ be vectors of non-negative integers. Then the $n_1 n_2\cdots n_m$ eigenvalues of $\sum\limits_{p=1}^{q} A^{k_{p1}}_1\otimes \cdots \otimes A^{k_{pm}}_m$ are $\xi_{j_1 j_2 \cdots j_m}=\sum\limits_{p=1}^{q} \xi^{k_{p1}}_{1j_1} \cdots \xi^{k_{pm}}_{mj_m}. $
 \\
For NEPS in graphs see \cite{cds,s}. The following definition extends this concept to sidigraphs.\\

\noindent {\bf Definition 4.2.} Let $\mathscr{B}$ be a set of binary $n-$tuples called basis for the product such that for every $i=1,2,\cdots,n$ there exists $\beta \in \mathscr{B} $ with $\beta_i=1$. The non-complete extended $p-$sum (or simply called NEPS) of sidigraphs $S_1,S_2,\cdots, S_n$ with basis $\bf \it B $ denoted by  NEPS$(S_1,S_2,\cdots,S_n;\mathscr{B})$ is a sidigraph with vertex set $V(S_1)\times V(S_2)\times \cdots\times V(S_n)$. There is an arc from vertex $(u_1,u_2,\cdots,u_n)$ to $(v_1,v_2,\cdots,v_n)$ if and only if there exists $(\beta_1,\beta_2,\cdots,\beta_n)\in \mathscr{B} $ such that $(u_i,v_i)\in A(S_i)$ whenever $\beta_i=1$ and $u_i=v_i$ whenever $\beta_i=0$. The sign of the arc is given by\\
$$\sigma((u_1,u_2,\cdots,u_n),(v_1,v_2,\cdots,v_n))=\prod\limits_{i=1}^{n} \sigma_i (u_i,v_i)^{\beta_i}=\prod\limits_{i:\beta_i=1}\sigma_i(u_i,v_i).$$
\indent Assume that the basis $ \mathscr{B} $ has $ r\ge 1 $ elements, i.e., $ \beta=\{\beta_1,\beta_2,\cdots,\beta_r\} \subseteq \{0,1\}^n\backslash \{(0,0,\cdots,0)\} $, we define  NEPS$(S_1,S_2,\cdots,S_n;\mathscr{B})=\bigcup\limits_{\beta\in \mathscr{B}}$ NEPS$(S_1,S_2,\cdots,S_n ;\beta).$\\

\noindent {\bf Example 4.3.}
The Kronecker product $S_1\otimes S_2 \otimes\cdots\otimes S_n $ of sidigraphs $S_1,S_2,\cdots,S_n$ is the NEPS of these sidigraphs with basis $\beta=\{(1,1,\cdots,1)\}$; the cartesian product $S_1\times S_2\times \cdots\times S_n$ is NEPS with basis $\{e_i\}, i=1,2,\cdots,n$, where $e_i$ is $n-$tuple with $1$ at $i$th position and $0$ otherwise.\\

The following result shows that two different basis give disjoint arc sets. The proof is straightforward.\\

\noindent {\bf Lemma 4.4.}
Let $S=$ NEPS $(S_1,S_2,\cdots,S_n;\beta)$ and $ S' =$ NEPS $ (S_1,S_2,\cdots,S_n;\beta')$, $\beta\ne \beta' $. Then $ \mathscr{A}(S)\cap \mathscr{A}(S')=\phi.$\\

The following result gives adjacency matrix and spectra of NEPS in terms of the constituent factor sidigraphs.\\

\noindent {\bf Theorem 4.5.}
If $S=$NEPS$(S_1,S_2,\cdots,S_n;\mathscr{B})$, then the ajacency matrix is given by\\
$$A(S)=\sum\limits_{\beta \in \mathscr{B}}A^{\beta_1}_1\otimes \cdots \otimes A^{\beta_n}_n,$$
and eigenvalues are given by  $z_{j_1 j_2 \cdots j_n}=\sum\limits_{\beta\in \mathscr{B}} z^{\beta_1}_{1j_1} \cdots z^{\beta_n}_{nj_n}$, where $1\le j_i\le n_i, i=1,2,\cdots,n.$\\
\noindent {\bf Proof.}
Let $ u=(u_{1 j_1},u_{2 j_2},\cdots,u_{n j_n})$ and $ v=(v_{1 j_1},v_{2 j_2},\cdots,v_{n j_n})$, where $1\le j_i \le |V(S_i)|$, $i=1,2,\cdots,n$ be any two vertices of $S$. Then
\begin{align*}
[A(S)]_{uv}&=\sum\limits_{\beta\in \mathscr{B}}(A^{\beta_1}_1)_{u_{1 j_1}v_{1 k_1}}(A^{\beta_2}_2)_{u_{2 j_2}v_{2 k_2}}\cdots (A^{\beta_n}_n)_{u_{n j_n}v_{n k_n}}
\\&=\sigma_1(u_{1j_1},v_{1k_1})^{\beta_1}\sigma_2(u_{2j_2},v_{2k_2})^{\beta_2}\cdots \sigma_n(u_{nj_n},v_{nk_n})^{\beta_n}
\\&=a^{\beta_1}_{1;j_1 k_1}a^{\beta_2}_{2;j_2 k_2}\cdots a^{\beta_n}_{n;j_n k_n}=[\sum\limits_{\beta \in \mathscr{B}}A^{\beta_1}_1\otimes \cdots \otimes A^{\beta_n}_n]_{uv}.
\end{align*}
The second part of the result follows by Lemma $4.1$. \qed \\

We note two special cases of Theorem $4.5$.\\
$(I)$ The Kronecker product $ S_1\otimes S_2\otimes \cdots \otimes S_n $ has eigenvalues $ z_{j_1 j_2 \cdots j_n}=z_{1j_1}z_{2j_2}\cdots z_{nj_n}$, $ 1\le j_i \le |V(S_i)|, i=1,2,\cdots,n $.\\
$(II)$ The cartesian product $ S_1\times S_2\times \cdots \times S_n $ has eigenvalues $ z_{j_1 j_2 \cdots j_n}=z_{1j_1}+ z_{2j_2}+\cdots + z_{nj_n}$, $ 1\le j_i \le |V(S_i)|, i=1,2,\cdots,n.$ \\

In \cite{gh2} the authors considered the problem of balance in NEPS of sigraphs. It is natural to consider the problem of  cycle balance for sidigraphs, in view of Theorem $1.3$, the cartesian product $S_1\times S_2\times \cdots \times S_n $ of sidigraphs $S_1,S_2,\cdots,S_n$, is  cycle balanced if and only if $S_i$, $i=1,2,\cdots,n$ is  cycle balanced. The next result gives sufficient but not necessary condition for cycle balance of NEPS and the proof follows on same lines as that in undirected case.\\

\noindent {\bf Theorem 4.6.} NEPS$(S_1,S_2,\cdots,S_n;\mathscr{B})$ is balanced if $S_1, S_2, \cdots, S_n$ are cycle balanced. \\

\noindent {\bf Remark 4.7} $(I)$ Theorem $4.6$ does not have a general converse. A counter example is $S=$ NEPS$(-C_3,- \overleftrightarrow{K_2},\{(1,1)\})$, where $-C_3$ denotes all negative directed cycle of order $3$ and $- \overleftrightarrow{K_2}$ is symmetric sidigraph of order $2$ with both arcs negative. $S$ is all positive, and hence cycle balanced. However $-C_3$ is non cycle balanced.\\
$(II)$ In view of Theorem $1.3$, the converse of Theorem $4.6$ is always true if basis $\mathscr{B}=\{e_i\}, i=1,2,\cdots,n$.\\

\noindent{\bf Theorem 4.8.} The following statements are equivalent about cartesian product $S=S_1\times S_2\times\cdots\times S_n$.\\
$(I)$ $S$ is cycle balanced.\\
$(II)$ All of $S_1,S_2,\cdots,S_n$ are cycle balanced.\\
$(III)$ $S$ and $S^u$ are cospectral.\\
\noindent{\bf Proof.} Theorem $1.3$ implies equivalence of $(I)$ and $(III)$ and Theorem $4.6$ implies equivalence of $(I)$ and $(II)$. \qed \\

\section{Upper bounds for the energy of sidigraphs}

Let $S$ be a sidigraph of order $n$ with adjacency matrix $A(S)=(a_{ij})$. The powers of $A(S)$ count the number of walks in signed manner. Let $w^+_{ij}(l)$ and $w^-_{ij}(l)$ respectively denote the number of positive and negative walks of length $l$ from $v_i$ to $v_j$. The following result relates integral powers of adjacency matrix with the number of positive and negative walks.\\

\noindent {\bf Theorem 5.1.} If $A$ is a adjacency matrix of a sidigraph on $n$ vertices, then $[A^l]_{ij}= w^+_{ij}(l)-w^-_{ij}(l)$. \\
\noindent {\bf Proof.} We prove the result by induction on $l$. For $l=1$, result is vacuously true. For $n=2$, let $n^+_{ij}$ denote the number of positive neighbours of distinct vertices $v_i$ and $v_j$, $n^-_{ij}$ the number of their common negative neighbours and $n^{\pm}_{ij}$ the number of neighbours that are positive to one vertex and negative to other. The $(i,i)$ entry of $A^2$ equals $w^+_{ii}(2)-w^+_{ii}(2)$. For $(i,j)$, $ i\neq j$,  $n^+_{ij}+ n^-_{ij}=w^+_{ij}(2)$ and  $n^{\pm}_{ij}=w^-_{ij}(2)$, so that $(i,j)$th entry$= w^+_{ij}(2)-w^-_{ij}(2)$. Now assume the result to be true for $l=m$.\\
We have, $[A^{m+1}]_{ij}=[A^m A]_{ij}=\sum\limits_{k=1}^{n}[A^m]_{ik}[A]_{kj}=w^+_{ij}(m+1)-w^-_{ij}(m+1)$,
by induction hypothesis. Therefore the result follows.  \qed \\
 
In sidigraph $S$, let $c^+_m$ denote number of positive closed walks of length $m$ and $c^-_m$ the number of negative closed walks of length $m$. In view of the fact that sum of eigenvalues of a matrix equals to its trace, we have following Corollary.\\

 \noindent {\bf Corollary 5.2.} If $z_1, z_2,\cdots,z_n$ are the eigenvalues of a sidigraph $S$, then $\sum\limits_{j=1}^{n}z^m_j=c^+_m-c^-_m$. \\

 \noindent {\bf Lemma 5.3.} Let $S$ be a sidigraph having $n$ vertices and $a$ arcs and let $z_1,z_2,\cdots,z_n$ be its eigenvalues. Then\\
 $(I)$ $\sum\limits_{j=1}^{n}(\Re z_j)^2-\sum\limits_{j=1}^{n}(\Im z_j)^2=c^+_2-c^-_2$, $(II)$ $\sum\limits_{j=1}^{n}(\Re z_j)^2+\sum\limits_{j=1}^{n}(\Im z_j)^2\le a=a^++a^-$.\\
 \noindent {\bf Proof.} By Corollary $5.2$ , we have\\
 $$ c^+_2-c^-_2=\sum\limits_{j=1}^{n} z_j^2=\sum\limits_{j=1}^{n}(\Re z_j)^2-\sum\limits_{j=1}^{n}(\Im z_j)^2+2\iota \sum\limits_{j=1}^{n}\Re z_j \Im z_j.$$\\
 Equating real and imaginary parts proves $(I)$.\\
 \indent By Schur's unitary triangularization, there exists a unitary matrix $U$ such that adjacency matrix $A$ of sidigraph $S$ is unitarily similar to an upper triangular matrix  $T=(t_{jk})$ with $t_{jj}=z_j$ for each $j=1,2,\cdots,n$. Then $\sum\limits _{j,k=1}^{n}|a_{jk}|^2=\sum\limits_{j,k=1}^{n}|t_{jk}|^2$. As $A$ is $(-1,0,1)-$matrix, we have\\
 $$a=\sum\limits_{j,k=1}^{n}|\sigma(v_j,v_k)|=\sum\limits_{j,k=1}^{n}|a_{jk}|=\sum\limits _{j,k=1}^{n}|a_{jk}|^2=\sum\limits_{j,k=1}^{n}|t_{jk}|^2\ge \sum\limits_{j=1}^{n}|t_{jj}|^2$$\\
  $$~~=\sum\limits_{j=1}^{n}|z_j|^2=\sum\limits_{j=1}^{n}\Re z_j^2+\sum\limits_{j=1}^{n}\Im z_j^2.$$\\
  thereby proving $(II)$. \qed \\

  \noindent {\bf Theorem 5.4.} Let $S$ be a sidigraph with $n$ verices and $a=a^+ +a^-$ arcs, and let $z_1,z_2,\cdots,z_n$ be its eigenvalues. Then $ E(S)\le \sqrt{\frac{1}{2}n(a+c^+_2-c^-_2)}.$
  \noindent {\bf Proof.} Subtracting part $(I)$ of Lemma $5.3$ from $(II)$, we see $\sum\limits_{j=1}^{n}(\Im z_j)^2\le \frac{1}{2}(a-(c^+_2-c^-_2))$. Applying Cauchy-Schwartz inequality to vectors $(|\Re z_1|,|\Re z_2|,\cdots,|\Re z_n|)$ and $(1,1,\cdots,1)$, we have\\
  $$E(S)=\sum\limits_{j=1}^{n}|\Re z_j|\le \sqrt{n} \sqrt{\sum\limits_{j=1}^{n}(\Re z_j)^2}=\sqrt{n}\sqrt{(c^+_2-c^-_2)+\sum\limits_{j=1}^{n}(\Im z_j)^2}$$\\
  $$\le \sqrt{n}\sqrt{(c^+_2-c^-_2)+\frac{1}{2}(a-(c^+_2-c^-_2)}=\sqrt{\frac{1}{2}n(a+c^+_2-c^-_2)}.$$ \qed \\

  \noindent {\bf Remark 5.5.} $(I)$ The upper bound in Theorem $5.4$ is attained by sidigraphs $S_1=(\frac{n}{2} \overleftrightarrow{K_2},+)$, $S_2=(\frac{n}{2} \overleftrightarrow{K_2},-)$, and skew-smmetric sidigraph of order $n$. Note that eigenvalues of $S_1$ and $S_2$ are $-1,+1$ each repeated $\frac{n}{2}$ times and eigenvalues of skew symmetric sidigraph of order $n$ are $0^{n-2},~\pm\iota\sqrt{n-1}$.\\
  $(II)$. Above result extends McClleland's inequality for sigraphs \cite{gh3} which states that $E(\Sigma)\le \sqrt{2pq}$, holds for every sigraph with $p$ vertices and $q$ edges. Let $\overleftrightarrow {\Sigma} $ be the symmetric sidigraph  of sigraph $\Sigma$, then in $\overleftrightarrow {\Sigma} $, $a=2q=c^+_2 = c^+_2-c^-_2$,. By Theorem $5.4$, $E(\Sigma)=E(\overleftrightarrow {\Sigma})\le \sqrt{\frac{1}{2}p(2q+2q)}=\sqrt{2pq}.$ \\
  
  The following result gives sharp upper bound of energy of sidigraphs in terms of number of arcs. The proof is same as in unsigned case.\\

  \noindent {\bf Theorem 5.6.} Let $S$ be a sidigraph with $a$ arcs. Then $E(S)\le a $ with equality if and only if $S=(\frac{a}{2}\overleftrightarrow {K_2},+)$ or $S=(\frac{a}{2}\overleftrightarrow {K_2},-)$  plus some isolated vertices.\\

  \noindent {\bf Remark 5.7.} Theorem $5.6$ extends the result for sigraphs \cite{gh3}, which states that $E(\Sigma)\le 2q$ for every sigraph with $q$ edges with equality if and only if $\Sigma= (\frac{q}{2}K_2,+)$ or  $\Sigma= (\frac{q}{2}K_2,-)$ plus some isolated vertices.

\section{Equienergetic sidigraphs}

Two sidigraphs are said to be isomorphic if their underlying digraphs are isomorphic such that signs are preserved. Any two isomorphic sidigraphs are obviously cospectral. There exist nonisomorphic sidigraphs which are cospectral, e.g., consider the sidigraphs $S_1$ and $S_2$ shown in Figure $2$. Clearly, $S_1$ and $S_2$ are nonisomorphic. $Spec~ S_1=\{0^{(5)}\}=Spec~S_2$.\\

\unitlength 1mm 
\linethickness{0.4pt}
\ifx\plotpoint\undefined\newsavebox{\plotpoint}\fi 
\begin{picture}(106.75,37.75)(0,0)
\put(19.75,18.25){\circle*{1}}
\put(19.75,18.25){\circle*{1}}
\put(19.5,33.25){\circle*{1.12}}
\put(19.5,33.25){\circle*{1.12}}
\put(33.25,26){\circle*{1}}
\put(33.25,26){\circle*{1}}
\put(49.75,32.75){\circle*{1.12}}
\put(49.75,32.75){\circle*{1.12}}
\put(48.5,16.75){\circle*{1.41}}
\put(48.5,16.75){\circle*{1.41}}
\put(82.25,32){\circle*{1}}
\put(82.25,32){\circle*{1}}
\put(81.5,17.75){\circle*{1}}
\put(81.5,17.75){\circle*{1}}
\put(98.5,31.5){\circle*{1}}
\put(98.5,31.5){\circle*{1}}
\put(26,29.63){\vector(2,-1){.07}}\multiput(19.25,33.25)(.0627907,-.03372093){215}{\line(1,0){.0627907}}
\put(26,29.63){\vector(2,-1){.07}}\multiput(19.25,33.25)(.0627907,-.03372093){215}{\line(1,0){.0627907}}
\put(41.5,29.5){\vector(2,1){.07}}\multiput(33.25,26)(.07932692,.03365385){208}{\line(1,0){.07932692}}
\put(41.5,29.5){\vector(2,1){.07}}\multiput(33.25,26)(.07932692,.03365385){208}{\line(1,0){.07932692}}
\put(19.5,26.13){\vector(0,1){.07}}\put(19.5,18.25){\line(0,1){15.75}}
\put(19.5,26.13){\vector(0,1){.07}}\put(19.5,18.25){\line(0,1){15.75}}
\put(89.88,32){\vector(1,0){.07}}\put(82,32){\line(1,0){15.75}}
\put(89.88,32){\vector(1,0){.07}}\put(82,32){\line(1,0){15.75}}
\put(98.25,25.13){\vector(0,-1){.07}}\put(98.25,31.5){\line(0,-1){12.75}}
\put(98.25,25.13){\vector(0,-1){.07}}\put(98.25,31.5){\line(0,-1){12.75}}
\put(102.5,34.5){\vector(3,2){.07}}\multiput(98.25,31.5)(.04775281,.03370787){178}{\line(1,0){.04775281}}
\put(102.5,34.5){\vector(3,2){.07}}\multiput(98.25,31.5)(.04775281,.03370787){178}{\line(1,0){.04775281}}
\put(26.38,22){\vector(-3,-2){.07}}\multiput(32.93,25.93)(-.7794,-.4706){18}{{\rule{.4pt}{.4pt}}}
\put(26.38,22){\vector(-3,-2){.07}}\multiput(32.93,25.93)(-.7794,-.4706){18}{{\rule{.4pt}{.4pt}}}
\put(48.88,24.75){\vector(0,-1){.07}}\multiput(49.43,32.93)(-.0735,-.9706){18}{{\rule{.4pt}{.4pt}}}
\put(48.88,24.75){\vector(0,-1){.07}}\multiput(49.43,32.93)(-.0735,-.9706){18}{{\rule{.4pt}{.4pt}}}
\put(40.63,21.25){\vector(-3,2){.07}}\multiput(48.18,16.43)(-.8026,.5){20}{{\rule{.4pt}{.4pt}}}
\put(40.63,21.25){\vector(-3,2){.07}}\multiput(48.18,16.43)(-.8026,.5){20}{{\rule{.4pt}{.4pt}}}
\put(81.75,25.13){\vector(0,-1){.07}}\multiput(81.93,32.18)(-.0333,-.95){16}{{\rule{.4pt}{.4pt}}}
\put(81.75,25.13){\vector(0,-1){.07}}\multiput(81.93,32.18)(-.0333,-.95){16}{{\rule{.4pt}{.4pt}}}
\put(106.25,37.25){\circle*{1}}
\put(106.25,37.25){\circle*{1}}
\put(31.5,14){\makebox(0,0)[cc]{$S_1$}}

\put(60,2.5){\makebox(0,0)[cc]{Figure 2}}
\put(60,2.5){\makebox(0,0)[cc]{Figure 2}}
\put(98,18.25){\circle*{1}}
\put(89.63,18.38){\vector(-1,0){.07}}\multiput(98,18.5)(-2.09375,-.03125){8}{\line(-1,0){2.09375}}
\put(87.25,11.5){\makebox(0,0)[cc]{$S_2$}}
\end{picture}

Two nonisomorphic sidigraphs $S_1$ and $S_2$ of same order are said to be equienergetic if $E(S_1)=E(S_2)$. In \cite{r} Rada proved the existence pairs of non-symmetric equienergetic digraphs on $n\ge 3$ vertices.  Cospectral sidigraphs are obviously equienergetic, therefore problem of equienergetic sidigraphs reduces to problem of construction of noncospectral pairs of equinergetic sidigraphs such that for every pair not both sidigraphs are cycle balanced. The problem of construction of equienergetic sigraphs is an open problem \cite{gh1} and for equienergetic graphs see \cite{abc,rw}.\\

\noindent {\bf Theorem 6.1.} Let $S$ be a sidigraph of order $n$ having eigenvalues $z_1,z_2,\cdots,z_n$ such that $|\Re z_j|\le 1$ for every $j=1,2,\cdots,n$. Then $E(S\times \overleftrightarrow{K_2})=2n$.\\
\noindent {\bf Proof.}
Let $z_1,z_2,\cdots,z_t$ be eigenvalues with nonnegative real part and $z_{t+1},\cdots,z_n$ be those with negative real part. Eigenvalues of cartesian product $S\times  \overleftrightarrow{k_2}$ are $z_1\pm 1,z_2 \pm 1,\cdots,z_t\pm 1,z_{t+1}\pm 1,\cdots,z_n\pm 1$. Therefore\\
$$ E(S \times \overleftrightarrow {K_2})=\sum\limits_{j=1}^{t}(|\Re z_j+1|+|\Re z_j-1|)+\sum\limits_{j=t+1}^{n}(|\Re z_j+1|+|\Re z_j-1|).$$
\indent As $|\Re z_j|\le 1$ for all $i=1,2,\cdots,n$, it follows that
$$E(S \times \overleftrightarrow {K_2})=\sum\limits_{j=1}^{t}(\Re z_j+1+1-\Re z_j)+\sum\limits_{j=t+1}^{n}(1-\Re z_j+\Re z_j+1)=2t+2(n-t)=2n.$$ \qed \\

\noindent {\bf Corollary 6.2.} For $ n\ge 2$, $E({\bf C}_n\times  \overleftrightarrow{K_2})=E(C_n\times \overleftrightarrow{K_2} )=2n$. Moreover, ${\bf C}_n\times  \overleftrightarrow{K_2}$ and $C_n\times \overleftrightarrow{K_2}$ are noncospectral sidigraphs with $2n$ vertices.\\
\noindent {\bf Proof.}
We know that eigenvalues of ${\bf C}_n$ are $e^{\frac{\iota (2j+1)\pi}{n}},~ j=0,1,\cdots,n-1$ and those of $C_n$ are $e^{\frac{2\iota j\pi}{n}},~ j=0,1,\cdots,n-1$. Clearly, eigenvalues of ${\bf C}_n$ and ${C_n}$ meet requirement of Theorem $6.1$, so $E({\bf C}_n\times  \overleftrightarrow{K_2})=E(C_n\times \overleftrightarrow{K_2} )=2n$. Moreover, $2\notin spec({\bf C}_n \times \overleftrightarrow{K_2})$ but $2\in spec(C_n\times \overleftrightarrow{K_2})$   implying that  ${\bf C}_n\times  \overleftrightarrow{K_2}$ and $C_n\times \overleftrightarrow{K_2}$ are non cospectral. Number of vertices in both sidigraphs is $2n$ follows by definition of cartesian product. In view of Theorem $1.3$,  ${\bf C}_n\times  \overleftrightarrow{K_2}$ is non cycle balanced whereas $C_n\times \overleftrightarrow{K_2}$ is cycle balanced.   \qed\\

\noindent {\bf Example 6.3.} For each odd $n$, ${\bf C}_n$ and $C_n$ is a non co-spectral pair of equienergetic sidigraphs, because $Spec({\bf C}_n)=-Spec(C_n)$ and $1\notin spec({\bf C}_n)$ but $1\in spec(C_n)$.\\

From Corollary $6.2$ and Example $6.3$, we see that for each positive integer $n\ge 3$  there exits a pair of non co-spectral sidigraphs with one sidigraph cycle balanced and another non cycle balanced. We now construct pairs of non cospectral sidigraphs of order $2n,~n\ge 5$ with both constituents non cycle balanced. Let $P^l_n~(n\ge l+1)$ be a sidigraph obtained by identifying one pendant vertex of the path $ P_{n-l+1} $ with any vertex of $ {\bf C}_l $. Sign of non cyclic arcs is immaterial.\\

\noindent {\bf Corollary 6.4.}  For each $n\ge 5$, $P^3_n\times \overleftrightarrow{K_2}$ and $P^4_n\times \overleftrightarrow{K_2}$ is a pair of noncospectral equienergetic sidigraphs of order and energy equal to $2n$.\\
\noindent {\bf Proof.} Using the fact that $\phi_{P^l_n}(x)=x^{n-l} \phi_{{\bf C}_l}(x)$ and Theorem $6.1$, it follows that $E(P^3_n\times \overleftrightarrow{K_2})=E(P^4_n\times \overleftrightarrow{K_2})=2n$. Now $1$ is an eigenvalue of  $P^3_n\times \overleftrightarrow{K_2}$ with multiplicity $n-3$ but $1$ is eigenvalue of $P^4_n\times \overleftrightarrow{K_2}$ with multiplicity $n-4$, therefore these two sidigraphs are noncospectral. Order of both sidigraphs equals to $2n$ follows by definition of cartesian product. In view of Remark $4.7~(II)$, it follows that both $P^3_n\times \overleftrightarrow{K_2}$ and $P^4_n\times \overleftrightarrow{K_2}$ are non cycle balanced. \qed \\

\noindent {\bf Corollary 6.5.} If $S$ is a sidigraph on $n$ vertices having eigenvalues $z_1,z_2,\cdots,z_n$ such that $|\Im z_j|\le \frac{1}{\sqrt{n-1}}$ for all $j=1,2,\cdots,n $, then $E((S\otimes {\bf S}_m)\times \overleftrightarrow{K_2})=2nm=$ order of sidigraph $(S\otimes {\bf S}_m)\times \overleftrightarrow{K_2}$, $2\le m\le n$.\\
\noindent {\bf Proof.} Eigenvalues of $S$ are $z_1,z_2,\cdots,z_n$ with $|\Im z_j|\le \frac{1}{\sqrt{n-1}}$ and eigenvalues of ${\bf S}_m $ are $0^{m-2},\pm \iota \sqrt{m-1}$. Then eigenvalues of Kronecker product $S\otimes {\bf S}_m $ meet requirement of Theorem $6.1$, therefore $E((S\otimes {\bf S}_m)\times \overleftrightarrow{K_2})=2nm$. Order of sidigraph equals to $2nm$ follows from definition of Kronecker product and cartesian product. \qed \\

\noindent{\bf Theorem 6.6.} Let $A=(a_{ij})$ be a square matrix of order $n$ having integral entries and zero trace and let $z_1,z_2,\cdots,z_n$ be its eigenvalues. Put $\alpha =\sum\limits_{j=1}^{n}|\Re z_j|$, then $\alpha $ cannot be of the form $(I)~ (2^t s)^{\frac{1}{h}}$ with $h\ge 1,~0\le t <h$ and $s$ odd  $(II)~ (\frac{m}{n})^{\frac{1}{r}}$ where $\frac{m}{n}$ is non-integral rational.\\
\noindent {\bf Proof.} We note that $\alpha = 2\sum\limits_{\Re z_j\ge 0}z_j$. Put $z=\sum\limits_{\Re z_j\ge 0}z_j$, then $z$ being sum of algebraic integers is an algebraic integer.\\
$(I)$ Assume $\alpha=(2^t s)^{\frac{1}{h}}$, so that $2z=(2^t s)^{\frac{1}{h}}$. Simplifying gives $z^h=\frac{s}{2^l}$, where $l=h-t\ge 1$. As $s$ is odd, therefore we see $z^h$ is non-integral rational algebraic integer, a contradiction. \\
$(II)$ Proof is same as in part $(I)$. \qed \\

Bapat and Pati \cite{bp} characterized positive rationals which can be the energy of graph. Later Pirzada and Gutman \cite{pg} proved energy of graph cannot be square root of an odd integer. These results were generalized to digraphs by Pirzada, Mushtaq, Gutman and Rada \cite{pmgr}. Following result generalizes these results to sidigraphs.\\

\noindent {\bf Theorem 6.7.} Energy of a sidigraph cannot be of the form  $(I)~ (2^t s)^{\frac{1}{h}}$ with $h\ge 1,~0\le t <h$ and $s$ odd  $(II)~ (\frac{m}{n})^{\frac{1}{r}}$ where $\frac{m}{n}$ is non-integral rational number and $r\ge 1$.\\
\noindent {\bf Proof.} Let $S$ be a sidigraph with adjacency matrix $A(S)$. Apply Theorem $6.6$ to $A(S)$ and note that $E(S)=\sum\limits_{j=1}^{n}|\Re z_j|$, the result follows. \qed \\

  \section{Open problems}
  
  We conclude this paper with the following open problems.\\
  $(1)$ Characterize sidigraphs with energy equal to the number of vertices.\\
  $(2)$ Find an infinite family of noncospectral equienergetic sidigraphs on $n\ge 4$ vertices with both constituents  non cycle balanced.\\
  $(3)$ Determine bases other than $\{e_i\}$ for which converse of Theorem $4.7$ holds.\\

  {\bf Acknowledgements}  The second author thanks University Grants Commission, New Delhi, India for providing junior research fellowship.\\

\end{document}